\documentclass[11pt,a4paper]{article}
\usepackage[margin=25mm]{geometry}
\usepackage{amsmath,amssymb,amsthm,graphicx}
\usepackage{microtype}
\usepackage{tikz}
\usepackage[colorlinks=true,citecolor=blue,urlcolor=blue]{hyperref}
\newtheorem{theorem}{Theorem}
\newtheorem{lemma}{Lemma}
\newcommand{\Z}{\mathbb Z}
\newcommand{\R}{\mathbb R}
\newcommand{\D}{\Delta}
\title{Straight boundary segments in the limiting shapes of Abelian sandpiles}
\author{Nikita Kalinin}
\date{}
\begin{document}
\maketitle
\begin{abstract}
Consider the Abelian sandpile on $\Z^2$ with $n$ grains added at the
origin to a constant integer background $h\le2$. Let $K$ be any
Hausdorff subsequential limit of the toppled sets rescaled by $n^{-1/2}$.
We prove that at the outermost point of $K$ on each coordinate half-axis $K$ agrees
locally with a closed half-plane. Thus the boundary contains a straight
segment about each of these four points and therefore cannot be a circle.
\end{abstract}

\section{Introduction}\label{sec:intro}
Fix an integer $h\le2$. Put $h$ grains at every site of $\Z^2$ and add
$n$ grains at the origin. A site containing at least four grains may
\emph{topple} by sending one grain to each of its four nearest neighbors.
For negative $h$, each site begins with a deficit of $-h$ grains.
Topplings continue until every site has height at most three. At these
background heights stabilization takes finitely many topplings and the
result is independent of their order \cite{BTW,Dhar,FLP}.

\begin{figure}[htbp]
 \centering
 \includegraphics[width=\textwidth]{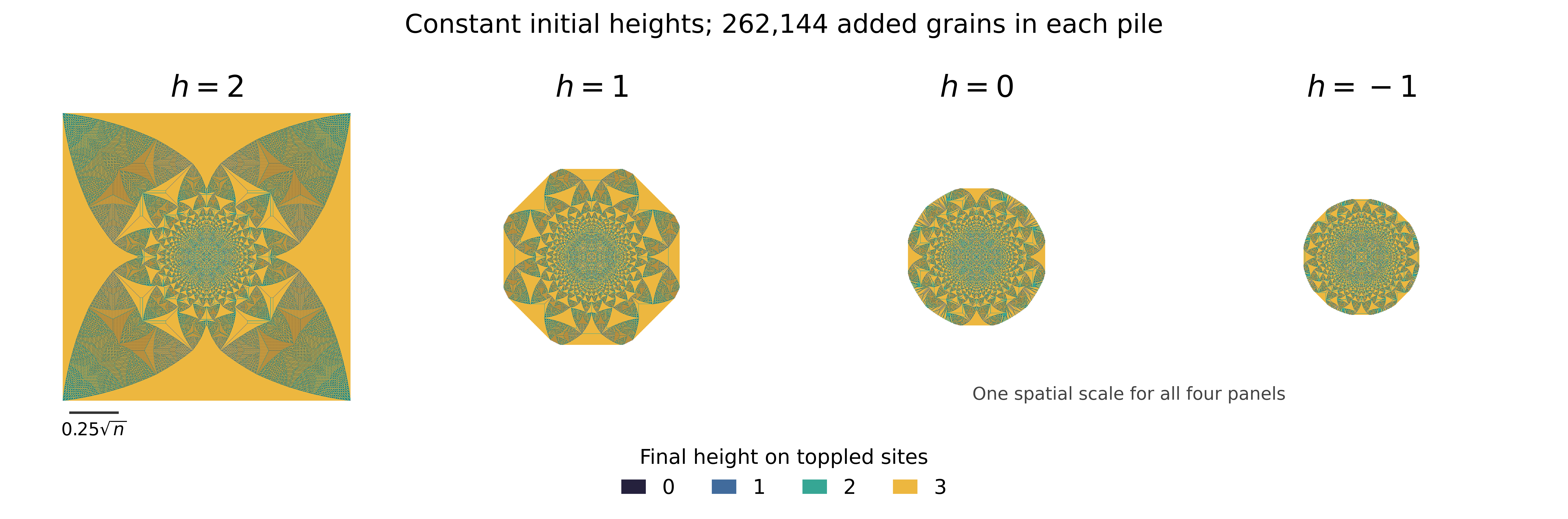}
 \caption{Exact relaxations of $n=262\,144$ added grains at the origin,
 with initial heights $h=2,1,0,-1$, from left to right. Colors give the
 final heights $0,1,2,3$ at sites that toppled. All other sites are white.
 Every panel has the same spatial scale.}
 \label{fig:backgrounds}
\end{figure}

The \emph{odometer} $u_n(x)$ is the number of topplings at $x$. We study the
shape of the toppled set $S_n=\{x:u_n(x)>0\}$ after dividing distances
by $\sqrt n$. Figure~\ref{fig:backgrounds} illustrates the difficulty of
reading the boundary from a finite configuration. At height two one can easily see a
square. At height zero, even the portions of the
boundary facing the coordinate axes look nearly round. So, this is the question addressed here: can a limiting boundary be an exact circle?

Fey-den Boer and Redig \cite[Theorem~4.1]{FeyRedig} proved that at
$h=2$ every nonempty toppled set is a lattice square. Holroyd, Levine,
M\'esz\'aros, Peres, Propp and Wilson \cite[Question~6.2]{HLMMPW}
asked whether the limiting aggregate is a regular $(12-4h)$-gon for
general $h\le2$. This would give an octagon at $h=1$ and a dodecagon
at $h=0$. In particular, they asked whether the apparent rounding at
the corners disappears in the limit. For the classical zero background,
Levine and Peres \cite{LP} discussed both a disk and a polygon as
possible shapes. Pegden and Smart \cite{PS} later singled out the
problem of ruling out a circular boundary.

Growth estimates determine the rescaling parameter.
For $h=0$, Le Borgne and Rossin \cite{LR} bounded $S_n$
from inside by a diamond and from outside by a square, both of diameter
of order $\sqrt n$. Levine and Peres \cite[Theorem~4.1]{LP} obtained an
inner ball bound for $h\le2$ and an outer ball bound for $h\le0$.
Fey, Levine and Peres completed the proof of the outer ball bound
for $h=1$ and obtained an outer square bound for $h=2$
\cite[Section~3.1 and Theorem~1.1]{FLP}. These estimates give the
$\sqrt n$ scale for each fixed $h$. The relative gap between the ball
bounds tends to zero as $h\to-\infty$. The spherical limit theorems in
\cite{LP} for the divisible sandpile and rotor-router aggregation
concern various different toppling rules.

For the zero background, Pegden and Smart \cite{PS} proved local
uniform convergence of the rescaled odometer away from the source and
weak-* convergence of the rescaled final density. Aleksanyan and
Shahgholian \cite[Theorem~4.7(b)]{AS} used directional monotonicity of
the limiting odometer to show that the boundary of its positivity set
is locally a Lipschitz graph. Bou-Rabee \cite[Theorem~9.2]{Bou} proved
local uniform odometer convergence away from the source for stationary
ergodic backgrounds bounded below and at most two, including every
constant height considered here. Levine, Pegden and Smart
\cite{Apollonian,LPS} constructed solutions of the sandpile equation
with Apollonian structure and classified the quadratic growth matrices
of integer superharmonic functions. Pegden and Smart \cite{Pattern}
proved pattern stability for certain classes of solutions. Passing
from odometer convergence to Hausdorff convergence of the toppled sets
would require further information wherever the limiting odometer vanishes.

Shmalo \cite[Theorem~1.1 and Corollary~1.2]{Shmalo} proved a
restriction on corners for every integer $h\le1$. If a Hausdorff limit
agrees near a point with a sector of opening $\alpha\in(0,\pi)$, then
$\alpha>\pi/2$ and $\sin\alpha\le 1/(3-h)$. Consequently a convex
polygonal limit would have strictly more than $12-4h$ vertices, ruling out
the proposed regular polygon. A disk has no corners so one cannot apply his methods.

Our result concerns the four points where the coordinate axes reach
the supporting lines of the limiting shape. At $h=2$ the conclusion follows
from the square theorem. The same conclusion holds at every fixed
integer $h\le2$, without a convexity or regularity assumption on the limit.

\begin{theorem}\label{thm:main}
Fix an integer $h\le2$ and let $K$ be any Hausdorff subsequential limit
of $n^{-1/2}S_n$. Write $R=\max\{x_1:x\in K\}$. Then $R>0$, and
there exists $\rho>0$ such that
\[
 K\cap B_\rho((R,0))
 =\{x\in B_\rho((R,0)):x_1\le R\}.
\]
The analogous statement holds at $(-R,0)$ and $(0,\pm R)$.
In particular, $\partial K$ contains a nontrivial straight segment
about each of these four points, and $K$ cannot be a Euclidean disk.
\end{theorem}

The starting point of the proof is a lower bound for the odometer along the horizontal axis. Writing \(d=R-x_1\), least action gives
\[
U(R-d,0)\geq\frac{3-h}{2}d^2.
\]
If this inequality holds in a small inward half-disk, then the positivity of the odometer implies that the limiting aggregate contains this half-disk. Its diameter is therefore a straight boundary segment.

To extend the inequality away from the axis, we subtract the quadratic term and study the resulting superharmonic function \(w\). A harmonic comparison first shows that the negative part of \(w\) is \(o(r^2)\) near the supporting point.

The remaining argument uses the integrality of the discrete odometer. We prove that certain quadratic polynomials cannot touch \(w\) from below. This is a particular instance of the Hessian constraints arising in the sandpile equation studied by Pegden, Smart and Levine, but the obstruction needed here follows from an elementary argument involving just two lattice steps. If \(w\) takes negative values arbitrarily close to the supporting point, we construct a quadratic polynomial violating this constraint.

\section{The odometer at the right-hand supporting point}
We use the unnormalized lattice Laplacian
$\D g(x)=\sum_{z\sim x}(g(z)-g(x))$.
\paragraph{A lower bound along the axis.}
Finite stabilization holds for every constant background $h\le2$. By least action, $u_n$ is the least
nonnegative integer-valued function $v$ satisfying
$n\delta_0+\D v\le3-h$ \cite[Lemmas~2.3--2.4]{FLP}.
In each comparison we lower some values of $u_n$ while keeping them
nonnegative and integral. If the result still satisfies the stability
inequality, least action forces equality with $u_n$ everywhere.

The first comparison shows that the largest value on each vertical
column is attained on the horizontal axis. The odometer inherits the
symmetries of the square from the initial configuration. 
Monotonicity of the odometer is classical\footnote{To keep the presentation self-contained, we may follow Shmalo's reflection argument for constant
backgrounds \cite[Lemma~2.1]{Shmalo}. For zero background this
monotonicity can be also found in Aleksanyan and Shahgholian \cite[Theorem~4.6]{AS}.
Reflect across $y=k+1/2$, where $k\ge0$
is an integer, and call the reflection $\sigma$. On $y\ge k+1$ replace
$u_n$ by $\min(u_n,u_n\circ\sigma)$. At each unchanged site its
Laplacian can only decrease. At a changed site $x$ the new value is
$u_n(\sigma x)$, and each neighboring value is at most its reflected
counterpart. Indeed, this follows from the minimum on the upper side.
If a neighbor lies across the reflecting line, it is $\sigma x$, whose
unchanged value $u_n(\sigma x)$ is strictly smaller than $u_n(x)$.
The new Laplacian at $x$ is therefore at most $\D u_n(\sigma x)\le3-h$.
All changed sites are away from the source. Least action forces
$u_n(i,k+1)\le u_n(i,k)$ on the first row above the line. Varying $k$
and using square symmetry gives monotonicity in each absolute coordinate.}.

For $n>3-h$ let $R_n$ be the largest integer $x$ for which $u_n(x,0)>0$.
Monotonicity gives $S_n\subset[-R_n,R_n]^2$. We claim that
\begin{equation}\label{axial}
 u_n(R_n-k,0)>\frac{3-h}{2}k(k+1),\qquad 0\le k\le R_n,\quad k\in\Z .
\end{equation}
To obtain the comparison function, look for a polynomial in $x$ with
lattice Laplacian $3-h$ and zeros in the columns $R_n$ and $R_n+1$.
This gives
\[
 V(x,y)=\frac{3-h}{2}(x-R_n)(x-R_n-1).
\]
Note that $V$ is nonnegative and integer-valued.
At $k=0$, \eqref{axial} follows from the definition of $R_n$.
If it failed at some $k>0$, monotonicity would give $u_n\le V$ on the
whole column $x=R_n-k$. Replace $u_n$ by $\min(u_n,V)$ strictly to
the right of this column. At each changed site the new values at all
four neighbors are bounded above by $V$, including any neighbor in
the comparison column. Hence the new Laplacian is at most $\D V=3-h$.
At unchanged sites, including the source, the Laplacian can only
decrease. This competitor vanishes at $(R_n,0)$, where $u_n$ is positive,
contradicting least action.

\paragraph{A local limit.}
Fix a subsequence along which $n^{-1/2}S_n\to K$ in Hausdorff distance.
The maximum of the first coordinate changes by at most that distance,
so $n^{-1/2}R_n\to R$. The points $n^{-1/2}(R_n,0)$ then give
$(R,0)\in K$, and symmetry gives $K\subset[-R,R]^2$.

Write $f_n(x)$ for the final height at $x$ minus $h$. A site that
topples has nonnegative height after its last toppling. At a site that
never topples the height stays at least $h$. Stability and the balance
of grains therefore give
\[
 \min(0,-h)\le f_n\le 3-h,\qquad \D u_n=f_n-n\delta_0.
\]
The support of $f_n$ is contained in $S_n$ and its nearest neighbors.
Summing the balance equation gives $\sum_x f_n(x)=n$, so
$n\le(3-h)(2R_n+3)^2$. In particular $R\ge1/(2\sqrt{3-h})>0$.
We can thus study the right edge in a ball that stays away from the source.

For completeness, we give a local compactness argument.  We use the potential
kernel as in Pegden and Smart \cite[Lemma~3.1]{PS} and Shmalo's
constant-background argument for $h\le1$ \cite[Proposition~3.1]{Shmalo}.
Let $A$ be the simple random walk potential kernel divided by four.
With our unnormalized Laplacian, $\D A=\delta_0$, and each coordinate
unit step $e$ satisfies $|A(x+e)-A(x)|\le C/(1+|x|)$
\cite[Proposition~4.4.2 and Corollary~4.4.5]{LL}. Since $u_n$ has
finite support, all sums in the following calculation are finite.
\[
 \begin{aligned}
 (A*\D u_n)(x)
 &=\sum_y A(x-y)\,\D u_n(y)\\
 &=\sum_y u_n(y)(\D A)(x-y)=u_n(x).
 \end{aligned}
\]
The balance equation now gives $u_n=A*f_n-nA$. For $x/\sqrt n$
within distance $R/2$ of $(R,0)$, the support of $f_n$ lies within
$C\sqrt n$ of $x$, while $|x|\ge R\sqrt n/2$. Taking one lattice
difference yields
\[
 |u_n(x+e)-u_n(x)|
 \le C\sum_{|z|\le C\sqrt n}\frac1{1+|z|}
       +\frac{Cn}{1+|x|}=O(\sqrt n).
\]
Indeed, the square shell $\|z\|_\infty=k$ contains $8k$ sites for
$k\ge1$, each with $|z|\ge k$. Thus for any $L\ge1$ the sum over
$|z|\le L$ is at most $1+\sum_{k=1}^{\lfloor L\rfloor}8k/(1+k)
\le1+8L$.

Place the value $u_n(x)/n$ at $x/\sqrt n$ and interpolate affinely on
the two right triangles of each mesh square. Along a horizontal or
vertical edge the value changes by $O(n^{-1/2})$ over a distance
$n^{-1/2}$, so both partial derivatives are uniformly bounded on each
triangle in the ball just considered. For large $n$, every triangle
meeting the closed ball of radius $R/4$ about $(R,0)$ lies in the
larger ball because its diameter is only $\sqrt{2/n}$. This gives a uniform
Lipschitz bound throughout the smaller ball.

At the fixed point $(R+R/8,0)$ the interpolant eventually vanishes.
Indeed, the point is at distance at least $R/8$ from $K$. By Hausdorff
convergence all vertices of a triangle containing it are eventually
outside $n^{-1/2}S_n$. The Lipschitz bound then bounds the interpolant
on the whole smaller ball. By Arzel\`a--Ascoli a further subsequence
converges uniformly there to a continuous function $U\ge0$.
The same argument at any point of this ball outside $K$ shows that $U$
vanishes there. In particular it vanishes on $x_1>R$, and by
continuity on $x_1=R$. By symmetry $U$ is even in the second coordinate. In particular, any
point of this neighborhood where $U>0$ must belong to $K$.

\paragraph{Subtracting the axial quadratic.}
\begin{figure}[htbp]
 \centering
 \begin{tikzpicture}[x=1cm,y=1cm,>=stealth,font=\small]
  \fill[black!4] (-3.9,-.9) rectangle (0,1.6);
  \draw[gray,dashed] (0,-.9) -- (0,1.85) node[above,text=black] {$x_1=R$};
  \draw[gray] (0,0) -- (2.0,0);
  \draw[->,thick] (0,0) -- (-3.65,0) node[above] {$d$};
  \draw[->,thick] (0,0) -- (0,1.42) node[right] {$y$};
  \draw[gray,dotted] (-2.45,0) -- (-2.45,.9) -- (0,.9);
  \fill (-2.45,.9) circle (1.7pt) node[above] {$(R-d,y)$};
  \fill (0,0) circle (2pt) node[below right] {$(R,0)$};
  \node at (-2.3,-.56) {$x_1<R$};
  \node at (1.2,.7) {$x_1>R$};
\end{tikzpicture}

 \caption{Coordinates at the right-hand supporting point. Increasing
 $d$ moves to the left. The shaded side is the half-plane containing
 $K$. The theorem shows that $K$ fills this side near $(R,0)$.}
 \label{fig:support}
\end{figure}
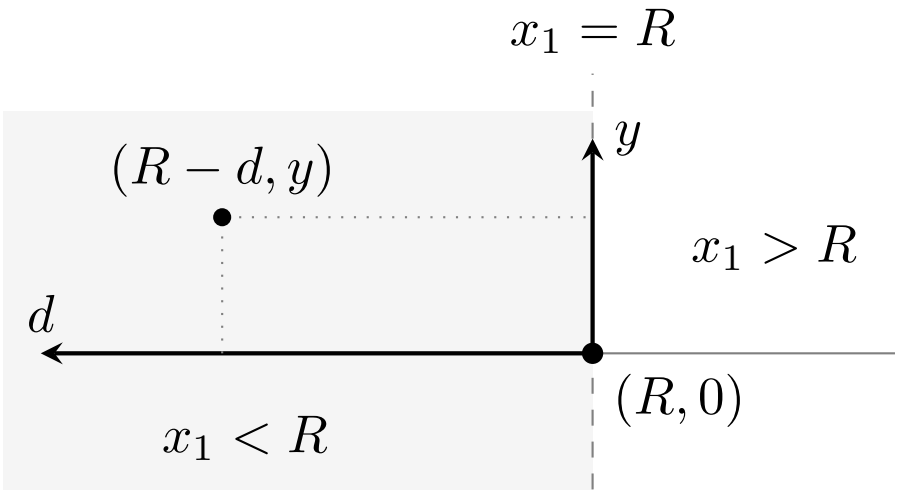
Use coordinates $(d,y)$ with $d$ pointing inward from $(R,0)$ (see Figure~\ref{fig:support}). Consider the error
\[
 w(d,y)=U(R-d,y)-\frac{3-h}{2}d^2 .
\]
It is important that the subtraction can be made on the lattice too.
At the mesh points $(d,y)=(k,j)/\sqrt n$ put
\[
 G_n(k,j)=u_n(R_n-k,j)-\frac{3-h}{2}k(k+1).
\]
The subtracted polynomial is integer-valued and has lattice Laplacian
$3-h$, so each $G_n$ is integer-valued and superharmonic. Moreover
$G_n(k,j)/n\to w(d,y)$ locally uniformly on the mesh, since
$R_n/\sqrt n\to R$ and $k(k+1)/n=d^2+d/\sqrt n$.

The limit $w$ is superharmonic as well. For a nonnegative smooth test
function $\phi$ compactly supported in this neighborhood, sum
$n^{-1}\phi(z/\sqrt n)\D G_n(z)$ over the lattice. The sum is
nonpositive. Summation by parts moves the Laplacian onto the sampled
test function, whose four differences sum to
$n^{-1}(\D\phi)(z/\sqrt n)+O(n^{-2})$, with $\D\phi$ the continuum
Laplacian. The resulting Riemann sums give $\int w\,\D\phi\le0$.

For small fixed $d>0$, take $k=\lfloor d\sqrt n\rfloor$ in
\eqref{axial} and pass to the limit. Together with the values on
the supporting line and the nonnegativity of $U$, this gives
\[
 w(d,0)\ge0,\qquad w(0,y)=0,\qquad
 w(d,y)\ge-\frac{3-h}{2}d^2.
\]
It remains to extend the first inequality into $d>0$.
By evenness in $y$ it suffices to work in the quadrant $d,y\ge0$.

\section{Extending the axial bound}
\paragraph{A first estimate on the error.}
We will prove
\begin{equation}\label{small}
 \max_{\substack{d,y\ge0\\d^2+y^2\le r^2}}(-w(d,y))_+=o(r^2)
 \qquad (r\downarrow0).
\end{equation}

Near the ray $d=0$, the lower bound $w\ge-\frac{3-h}{2}d^2$ already
suffices. Fix $0<\tau<1$. In the narrow region $d<\tau y$ of the
quarter disk we have $-w\le\frac{3-h}{2}\tau^2r^2$. The remaining
sector, between $y=0$ and $d=\tau y$, has the angle
$\theta_0=\arctan(1/\tau)<\pi/2$. Adding the harmonic function
$\frac{3-h}{2}\tau dy$ to $w$ gives a superharmonic function
nonnegative on both rays. On $y=0$ the added term vanishes, and on
$d=\tau y$ it exactly cancels $-\frac{3-h}{2}d^2$.

By the strict inequality $\theta_0<\pi/2$ the sector admits a positive homogeneous harmonic function of degree greater than two. Choose $2<p<\pi/\theta_0$.
In polar coordinates $(s,\theta)$ measured from the $d$-axis, take
\[
 H(s,\theta)=s^p\cos\bigl(p(\theta-\theta_0/2)\bigr),
\]
which is positive on the sector. A
sufficiently large multiple $C_\tau H$ therefore makes
$w+\frac{3-h}{2}\tau dy+C_\tau H$ nonnegative on a fixed small outer
arc. It is also nonnegative on the rays and at the vertex. The minimum
principle in the sector and $dy\le s^2/2$ give
\[
 -w\le\frac{3-h}{4}\tau s^2+C_\tau s^p .
\]
For fixed $\tau$, the last term is $o(r^2)$ uniformly for $s\le r$.
Combining the two sectors, dividing by $r^2$, and then letting $\tau$
decrease to zero proves \eqref{small}.

\paragraph{Quadratics touching the limit from below.}
Suppose a quadratic lies below $w$ on the boundary of a domain but
above it somewhere inside. Lowering the quadratic by a constant until
it lies below $w$ everywhere produces an interior point of equality.
We now show that, for the curvatures in Lemma~\ref{lem:two}, the
lattice rules out this possibility even in the limit. This comparison
is related to the Hessian constraints studied in \cite[Lemma~4.1]{PS},
\cite[Section~2.2 and Proposition~7]{Pattern}, and
\cite[Theorem~1.1]{LPS}.

\begin{lemma}\label{contact}
Let $D\subset\R^2$ be bounded, open and nonempty. Suppose that along
a sequence $n\to\infty$ the mesh values $G_n(z)/n$ at $z/\sqrt n$
of integer-valued superharmonic functions converge locally uniformly
to a continuous function $W$ in a neighborhood of $\overline D$.
If the Hessian of a quadratic polynomial $P$ satisfies the hypotheses
of Lemma~\ref{lem:two}, then
\[
 \min_{\overline D}(W-P)=\min_{\partial D}(W-P).
\]
In particular, $P\le W$ on $\partial D$ implies $P\le W$ in $D$.
\end{lemma}

\begin{proof}
Suppose the interior minimum is strictly smaller than the boundary
minimum. Choose a lattice point $z_n$ at which $G_n(x)-nP(x/\sqrt n)$
is smallest among the points with $x/\sqrt n\in\overline D$. Mesh points
approaching an interior minimizer of $W-P$ give rescaled values tending
to its minimum. A subsequence with $z_n/\sqrt n$ approaching the
boundary would give a strictly larger limiting value, by continuity and
local uniform convergence. Thus $z_n/\sqrt n$ stays a fixed positive
distance from $\partial D$. For large $n$ all points $(z_n+z)/\sqrt n$
with $|z|\le2$ lie in $D$.

On these sites the minimizing inequality and the exact quadratic
Taylor expansion, with $B=D^2P$, give
\[
 \begin{aligned}
 G_n(z_n+z)-G_n(z_n)
 &\ge n\left[P\left(\frac{z_n+z}{\sqrt n}\right)
             -P\left(\frac{z_n}{\sqrt n}\right)\right]\\
 &=\sqrt n\,\nabla P(z_n/\sqrt n)\cdot z+\tfrac12 z^TBz .
 \end{aligned}
\]
We remove the linear term up to an error of at most \(1/2\) in each coordinate. Subtracting the nearest
integer slope preserves integrality and the Laplacian. For this
fixed $n$, write $a=\sqrt n\,\nabla P(z_n/\sqrt n)$ and set
$m_i=\lfloor a_i+1/2\rfloor$, $b=a-m$. Then $|b_i|\le1/2$ and
\[
 g(z)=G_n(z_n+z)-G_n(z_n)-m\cdot z
 \ \ge\ \tfrac12 z^TBz+b\cdot z\qquad (|z|\le2).
\]
Here $g(0)=0$, the integer vector $m$ makes $g$ integer-valued, and
the affine function subtracted from $G_n$ has zero lattice Laplacian.
Thus $g$ is superharmonic and violates Lemma~\ref{lem:two} below.
\end{proof}

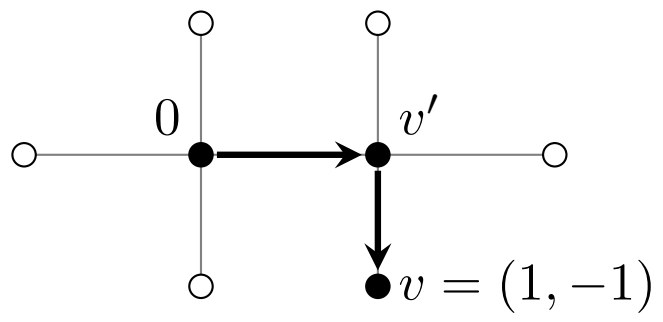
\begin{figure}[htbp]
 \centering
 \begin{tikzpicture}[x=1.17cm,y=.87cm,>=stealth,font=\small]
  \draw[gray] (-1,0) -- (2,0) (0,-1) -- (0,1) (1,-1) -- (1,1);
  \foreach \x/\y in {-1/0,0/1,0/-1,1/1,2/0}
    \draw[fill=white] (\x,\y) circle (2.2pt);
  \draw[->,very thick,shorten <=3pt,shorten >=3pt] (0,0) -- (1,0);
  \draw[->,very thick,shorten <=3pt,shorten >=3pt] (1,0) -- (1,-1);
  \fill (0,0) circle (2.4pt) node[above left] {$0$};
  \fill (1,0) circle (2.4pt) node[above right] {$v'$};
  \fill (1,-1) circle (2.4pt) node[right] {$v=(1,-1)$};
\end{tikzpicture}

 \caption{One choice of the two steps. The thin edges join $0$ and
 $v'$ to the neighbors checked in the proof. The arrows trace the
 propagation of zeros from $0$ to $v'$ and then to $v$.}
 \label{fig:steps}
\end{figure}

\begin{lemma}\label{lem:two}
Let $B$ be a symmetric $2\times2$ matrix with $B\succeq-I/8$ and
$v^TBv>0$ for $v=(1,-1)$. For any $b$ with $|b_i|\le1/2$, there
is no integer-valued function $g$ on the lattice ball $|z|\le2$
that is superharmonic at $|z|\le1$ and satisfies
\[
 g(0)=0,\qquad g(z)\ge q(z):=\tfrac12 z^TBz+b\cdot z
 \quad (|z|\le2).
\]
\end{lemma}

\begin{proof}
Changing the sign of $v$ preserves its quadratic value, so choose the
sign for which $b\cdot v\ge0$. The two coordinate steps from $0$ to
$v$ have $b$-increments summing to $b\cdot v$. Take a step with
nonnegative increment first, and call its endpoint $v'$
(Figure~\ref{fig:steps}). Every neighbor $z$ of $0$ or $v'$ has
$|z|\le2$. A single coordinate step lowers the linear term by at most
$1/2$, and that term is nonnegative at both centers. Consequently
\[
 g(z)\ge q(z)\ge-\frac{|z|^2}{16}-\frac12\ge-\frac34>-1
\]
at all their neighbors, where the quadratic bound uses $B\succeq-I/8$.
 Suppose an integer-valued
superharmonic function $g$ vanishes at $x$ and $g(z)>-1$ at its four
neighbors. Integrality gives $g(z)\ge0$ there, while superharmonicity
gives $\sum_{z\sim x}g(z)\le4g(x)=0$. This first forces $g(v')=0$ and then forces $g(v)=0$,
since $v$ is a neighbor of $v'$. Yet $q(v)>0$ by the choice of sign
and the strict positivity of $v^TBv$.
\end{proof}

\begin{figure}[htbp]
 \centering
 \begin{tikzpicture}[x=1.32cm,y=1.1cm,>=stealth,font=\small]
  \fill[black!4] (0,0) -- (3,0) -- (0,3) -- cycle;
  \draw[->] (0,0) -- (3.57,0) node[right] {$d$};
  \draw[->] (0,0) -- (0,3.45) node[above] {$y$};
  \draw[thick] (0,0) -- (3,0) -- node[midway,above right] {$d+y=2r$} (0,3) -- cycle;
  \draw[dashed] (0,1.5) -- (1.5,0);
  \node[anchor=east] at (0,1.5) {$r$};
  \node[anchor=east] at (0,3) {$2r$};
  \node[below] at (1.5,0) {$r$};
  \node[below] at (3,0) {$2r$};
  \node[below left] at (0,0) {$0$};
  \draw[gray,thin] (.92,.58) -- (1.08,.53);
  \node[anchor=west,inner sep=1pt] at (1.08,.53) {$d+y=r$};
  \fill (.44,1.06) circle (2pt) node[above right] {$(d_0,y_0)$};
\end{tikzpicture}
 \caption{The dashed segment contains the point where $w<0$.
 On the axes $w\ge0$. On the third edge the first error estimate
 bounds how far $w$ can fall below zero.}
 \label{fig:triangle}
\end{figure}
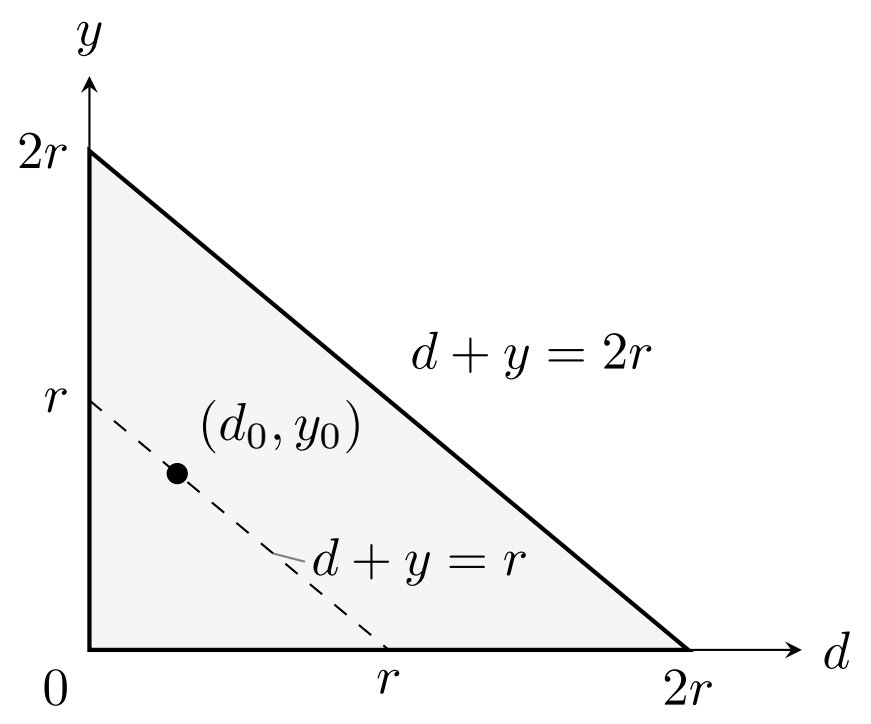

\paragraph{Why the boundary is straight, proof of Theorem~\ref{thm:main}.}
Suppose $w$ takes negative values arbitrarily close to $(0,0)$ in the
quadrant $d,y\ge0$. We know that both coordinates of any such point
are positive. Choose one, $(d_0,y_0)$, and put $r=d_0+y_0$. The line
$d+y=r$ passes through it and meets both axes, where we already
control $w$. Enclose this line segment in the triangle
$D=\{d,y>0:d+y<2r\}$ as in Figure~\ref{fig:triangle}.

We want the polynomial to vanish on the dashed segment, so that it
lies above $w$ at the chosen point, and to decrease as we move away
from that segment. Squaring the distance to the line, up to a constant
factor, suggests
\[
 P_0(d,y)=-\frac{(d+y-r)^2}{64}.
\]
This polynomial is nonpositive on the axes and equals $-r^2/64$ on
the third edge. That edge lies within distance $2r$ of the vertex
$(0,0)$. By \eqref{small}, the maximum of $(-w)_+$ is $o(r^2)$ there so we can choose the negative point sufficiently close to the origin that the
closed triangle lies in our coordinate neighborhood and $P_0\le w$
on the third edge.

Along $(1,-1)$, the direction of its line of zeros, $P_0$ has zero
curvature. Subtracting $\varepsilon dy$ gives positive curvature
in that direction, leaves both coordinate edges unchanged, and only
lowers the polynomial on the third edge. For $P=P_0-\varepsilon dy$
a direct differentiation gives
\[
 D^2P=
 -\frac1{32}\begin{pmatrix}1&1\\1&1\end{pmatrix}
 -\varepsilon\begin{pmatrix}0&1\\1&0\end{pmatrix}.
\]
The eigenvalues along $(1,-1)$ and $(1,1)$ are respectively
$\varepsilon$ and $-1/16-\varepsilon$. For \(0<\varepsilon<1/16\), the Hessian therefore satisfies the required lower bound. Choose
\[
 0<\varepsilon<\min\left\{\frac1{16},
                 \frac{-w(d_0,y_0)}{d_0y_0}\right\}.
\]
Then $D^2P\succeq-I/8$ and $P(d_0,y_0)=-\varepsilon d_0y_0>w(d_0,y_0)$.
Since $P\le w$ on all three edges, Lemma~\ref{contact} with $W=w$
gives the desired contradiction.

Thus $w\ge0$ in the part of a small disk where $d,y\ge0$.
Reflection in $y=0$ gives the same bound in the entire inward half-disk.
For some $\rho>0$ we have
\[
 U(R-d,y)\ge\frac{3-h}{2}d^2>0
 \qquad (d>0,\ d^2+y^2<\rho^2).
\]
Every point in this open half-disk belongs to $K$, because $U$
vanishes outside $K$. Every point of its diameter inside the ball
also belongs to $K$, by closedness, as a limit of points from the
inward side. On the other side $x_1>R$, which excludes membership
in $K$ by the definition of $R$. We have proved
\[
 K\cap B_\rho((R,0))=\{x\in B_\rho((R,0)):x_1\le R\}.
\]
Rotations by multiples of $\pi/2$ preserve each toppled set and hence
$K$. They give the other three straight segments in
Theorem~\ref{thm:main}. \hfill$\square$


\begin{thebibliography}{99}
\small
\setlength{\parsep}{0pt}
\setlength{\itemsep}{3pt}
\bibitem{AS} H. Aleksanyan, H. Shahgholian,
\emph{Discrete balayage and boundary sandpile},
J. Anal. Math. \textbf{138} (2019), 361--403.
\url{https://arxiv.org/abs/1607.01525}.

\bibitem{BTW} P. Bak, C. Tang, K. Wiesenfeld,
\emph{Self-organized criticality: an explanation of the $1/f$ noise},
Phys. Rev. Lett. \textbf{59} (1987), 381--384.
\url{https://doi.org/10.1103/PhysRevLett.59.381}.

\bibitem{Bou} A. Bou-Rabee,
\emph{Convergence of the random Abelian sandpile},
Ann. Probab. \textbf{49} (2021), 3168--3196.
\url{https://arxiv.org/abs/1909.07849}.

\bibitem{Dhar} D. Dhar,
\emph{Self-organized critical state of sandpile automaton models},
Phys. Rev. Lett. \textbf{64} (1990), 1613--1616.
\url{https://doi.org/10.1103/PhysRevLett.64.1613}.

\bibitem{FLP} A. Fey, L. Levine, Y. Peres,
\emph{Growth rates and explosions in sandpiles},
J. Stat. Phys. \textbf{138} (2010), 143--159.
\url{https://arxiv.org/abs/0901.3805}.

\bibitem{FeyRedig} A. Fey-den Boer, F. Redig,
\emph{Limiting shapes for deterministic centrally seeded growth models},
J. Stat. Phys. \textbf{130} (2008), 579--597.
\url{https://arxiv.org/abs/math/0702450}.

\bibitem{HLMMPW} A. E. Holroyd, L. Levine, K. M\'esz\'aros,
Y. Peres, J. Propp, D. B. Wilson,
\emph{Chip-firing and rotor-routing on directed graphs},
in \emph{In and Out of Equilibrium 2}, Progress in Probability
\textbf{60}, Birkh\"auser, 2008, 331--364.
\url{https://arxiv.org/abs/0801.3306}.

\bibitem{LL} G. F. Lawler and V. Limic,
\emph{Random Walk: A Modern Introduction},
Cambridge University Press, 2010.

\bibitem{LR} Y. Le Borgne and D. Rossin,
\emph{On the identity of the sandpile group},
Discrete Math. \textbf{256} (2002), 775--790.
\url{https://doi.org/10.1016/S0012-365X(02)00347-3}.

\bibitem{Apollonian} L. Levine, W. Pegden, C. K. Smart,
\emph{Apollonian structure in the Abelian sandpile},
Geom. Funct. Anal. \textbf{26} (2016), 306--336.
\url{https://arxiv.org/abs/1208.4839}.

\bibitem{LPS} L. Levine, W. Pegden, C. K. Smart,
\emph{The Apollonian structure of integer superharmonic matrices},
Ann. of Math. \textbf{186} (2017), 1--67. Corrected author version:
\url{https://www.math.cmu.edu/~wes/pdfs/integer-superharmonic.pdf}.

\bibitem{LP} L. Levine, Y. Peres,
\emph{Strong spherical asymptotics for rotor-router aggregation and
the divisible sandpile},
Potential Anal. \textbf{30} (2009), 1--27.
\url{https://arxiv.org/abs/0704.0688}.

\bibitem{PS} W. Pegden, C. K. Smart,
\emph{Convergence of the Abelian sandpile},
Duke Math. J. \textbf{162} (2013), 627--642.
\url{https://arxiv.org/abs/1105.0111}.

\bibitem{Pattern} W. Pegden, C. K. Smart,
\emph{Stability of patterns in the Abelian sandpile},
Ann. Henri Poincar\'e \textbf{21} (2020), 1383--1399.
\url{https://arxiv.org/abs/1708.09432}.

\bibitem{Shmalo} Y. Shmalo,
\emph{Corner angles of sandpile aggregates on a background of holes},
preprint dated 13 October 2026.
\href{https://github.com/yspennstate/sandpile-holes-corner-angles}{Paper and companion code}.

\end{thebibliography}
\end{document}